\documentclass[reqno, 10pt]{amsart}
\usepackage{color}

\usepackage{graphicx, enumerate}
\usepackage{amssymb, amsmath, amsthm}
\allowdisplaybreaks

\usepackage{hyperref}
\hypersetup{
    colorlinks=true, 
    linktoc=all,     
    linkcolor=blue,  
    citecolor=blue
}

\newtheorem{Theorem}{Theorem}[section]

\newtheorem{Definition}[Theorem]{Definition} 

\theoremstyle{definition}

\theoremstyle{remark}

\numberwithin{equation}{section}

\begin{document} 

\title[partitions as words]{Lognormal degree distribution in the partition graphs}
\author[H.~S.~Bal]{Hartosh Singh Bal
}
\address{The Caravan\\
Jhandewalan Extn., New Delhi 110001, India}
\email{hartoshbal@gmail.com}

\date{\today}


\maketitle

\begin{abstract}
We demonstrate a method for listing all ordinary partitions of $n$ as binary words of length $(n-1)$. The resulting family $\Pi_n$ imbued with the hamming distance yields subgraphs of the Hamming Graphs. The existence of a 2-Gray Code for ordinary partitions follows from the fact that the graph (with the $0^{n-1}$ partition omitted) is 2-connected. However, the graphs fail to be hamiltonian for ordinary partitions when $n>7$, ruling out the possibility of a Gray code for all such flip graphs. We further investigate the degree distribution of the graph for $n$, and provide computational evidence that this is a long-tailed lognormal distribution. This conjecture connects to a closely related, and much older, question of the distribution of the number of parts of a partition and the same evidence suggests that this distribution is also lognormal for large $n$. These methods extend to higher dimensional partitions of $n$ which can be then written as words of length $(n-1)$ on $d+1$ letters. The resulting graphs are connected, proving that d-dimensional partitions allow a 3-Gray code. 
\end{abstract}

\section{Partitions as binary words}\label{sec:intro}

The partition function $p(n)$ counts the number of ways in which we can write $n$ as a sum  $\sum_{i=1}^{n}\lambda_i = n$ where the $\lambda_i$ are non-increasing positive integers. We begin by listing out rules which uniquely describe each ordinary partition as a special category of binary words. \par

 \begin{Definition}\label{rule1}       Rules for binary partition strings: \\
        a) Consecutive $0$s can only terminate a string.\\
        b) Any string of $1$s is allowed and is termed a block.\\
        c) Given $A_1 = 1^{\lambda_1}$ and  $A_2=1^{\lambda_2}$, $A_10A_2$ is allowed if				  $\lambda_1\geq\lambda_2$.\\
        d)Given blocks $A_i$ and an allowed word $A_10A_20\cdots0A_n$, the word $A_10\cdots0A_n0A_{n+1}$ is allowed if $A_n0A_{n+1}$ is allowed.\\
        \end{Definition}
  Beginning with the empty word, we recursively generate allowed words, appending $0$ to all allowed words of length $n-1$, and appending $1$ to all allowed words of length $n-1$ where the rules allow. In the table, $p(n)$ denotes the number of partitions of $n$. 
$$
\begin{array}{| c | l | c |}\hline
n & \text{Allowed words of length  } n-1 & p(n) \\ \hline\hline
1 & \epsilon & 1\\ 
2 &  0, 1 & 2 \\ 
3 & 00, 10, 11 & 3 \\
4 & 000, 100, 110, 101, 111 & 5\\
5 & 0000, 1000, 1100, 1010,1110, 1101, 1111 & 7 \\
6 &  \begin{minipage}[t]{0.75\textwidth}
$00000,$ $10000,$ $11000,$ $10100,$ $11100,$ $11010,$ $11110,$ $10101,$\\
$11101,$ $11011,$ $11111$
\end{minipage} & 11\\
\hline
\end{array}
 $$
 
We term a word allowed by these rules a partition word. A generic partition word is of the form $1^{\lambda_1}01^{\lambda_2}\cdots1^{\lambda_i}01^{\lambda_{i+1}}\cdots1^{\lambda_n}0^m$ where ${\lambda_i}\geq{\lambda_{i+1}}$ for all $i$, $m$ is any arbitrary non-negative integer and the expression $a^0$ represents the empty word. \par
  
\begin{Theorem}\label{ordinary}  The rules listed in the introduction ~\eqref{rule1} give a bijection between all partitions of $n$ and allowed binary words of length $n-1$.
\end{Theorem}
\begin{proof}
Consider any partition $\sum_{i=1}^{n}\lambda_i = n$ where the $\lambda_i$ are non-increasing positive integers, and we assume without loss of generality that for some $1 \leq j \leq n$, all the $\lambda_k$ for $k \geq j$ are $1$. Then, the allowed partition word $1^{\lambda_1 - 1}01^{\lambda_2 - 1}0 \cdots 01^{\lambda_{j-1} - 1}0^{n-j}$ corresponds uniquely to the given partition.
\end{proof}\par
The number of $0$s in the partitions of $n$ correspond the number of plus signs in the partitions, given by OEIS A076276, see \cite{sloane}, while the number of $1$s correspond to the sum of all parts of a partition minus the number of parts, given by OEIS A196087, see ~\cite{sloane}. By our construction the sum of these two quantities adds up to $(n-1)p(n)$ for each $n$.\par
The basic idea here goes back to Macmahon ~\cite{macmahon1} who had proposed it in the context of compositions. In a recent paper Sills ~\cite{sills} has observed how this defines a unique word for each composition. Our rules give an effective way of determining when such a composition word is a partition. The role played by $0$ and $1$ in Sills' paper has  been interchanged here. \par

\section{Partition Subgraphs of the Hypercube}\label{graph}

We define the distance between any two strings as the number of single letter alterations needed to obtain one from the other. For example, the distance, termed the Hamming distance, between $101$ and $110$ is $2$. We consider the set of binary partition words for each $n$ as vertices of the graph $\Pi(1,n)$, where all points separated by a distance of $1$ are connected by edges.\\

\begin{Theorem}
$\Pi(1,n)$ is a connected bipartite subgraph of the hypercube $Q_{n-1}$. 
\end{Theorem}
\begin{proof}
The partition graph for $n=2$, consisting of words of length $1$, is a connected bipartite subgraph of $Q_1$. If we assume the result is true for $n-1$, i.e. words of length $n-2$, then the partition words of length $n-1$ are of two types, those obtained from partition words of length $n-2$ by appending $0$ to each word, and those obtained by appending $1$ where allowed. The words ending in $0$ clearly form a connected subgraph of $Q_{n-1}$. Moreover, any word ending in $1$ will also have an associated word ending in $0$ which only differs at the last place, and hence it will be connected to the subgraph by a vertex of the $Q_{n-1}$ hypercube. It is also clear that the graph is bipartite because the words with an odd(even) number of $1$s only connect the words with an even(odd) number of $1$s.
\end{proof}

\textbf{The graph of P(8), consisting of 22 vertices, is illustrated below.}

\resizebox{\textwidth}{!}{\includegraphics{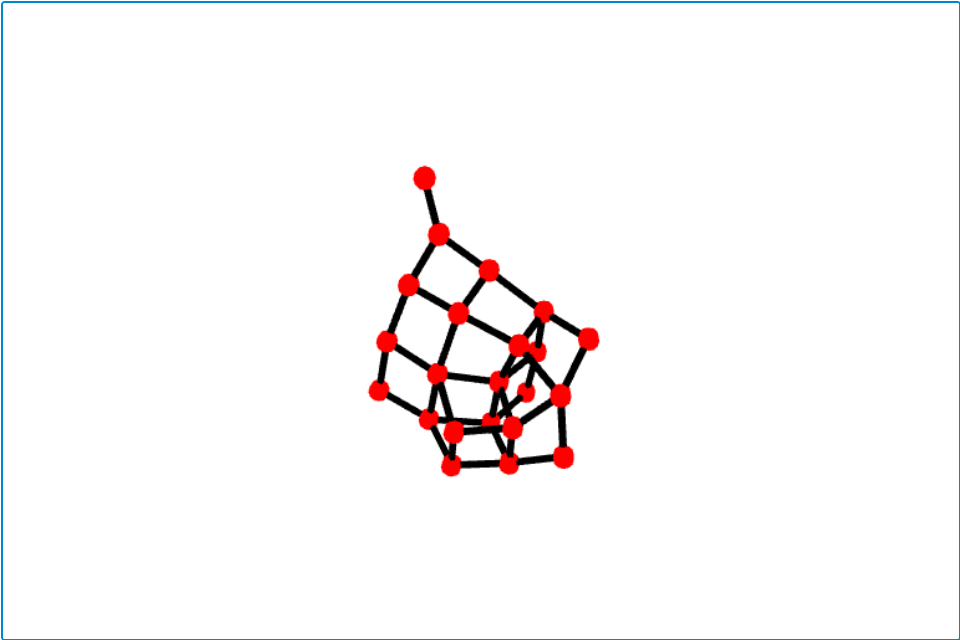}}

\begin{Theorem}
The diameter of $\Pi(1,n)$ is $n-1$.
\end{Theorem}
\begin {proof}
The distance between the vertices $A=0^{n-1}$ and $B=1^{n-1}$ is $n-1$. We note that any vertex has a connected path to $A$ through allowed partition words by moving right to left in the string and changing each $1$ to a $0$ or to $B$ by moving left to right in the string and changing each $0$ to a $1$. Let $x_1$ and $x_2$ be any two vertices of the graph. If the distance of $x_1$ from $A$ is $a_1$, then the distance of $x_1$ from $B$ is $n-1-a_1$, and one of the two terms is $\le$ $\lfloor{(n-1)/2}\rfloor$, similarly the distance of $x_2$ from $A$ is $b_1$, then the distance of $x_2$ from $B$ is $n-1-a_2$, and one of them is $\le$  $\lfloor{(n-1)/2}\rfloor$. By choosing the appropriate segment joining $x_1$ to $A$ or $B$, as the case may be, it follows that the shortest path between $x_1$ and $x_2$ via $A$ or $B$ is $\leq$ ${n-1}$.
\end{proof}

Given that the number of vertices $N$ of  $\Pi(1,n)$ is $\mathcal{O}(e^{\sqrt{n}})$, we get that asymptotically the diameter of $\Pi(1,n)$ is $\mathcal{O}(\log^{2}{N})$

\begin{Theorem}
The graph $\Pi_n / \lbrace 0^{n-1}\rbrace$ (i.e. with the $0$ vertex removed) is a 2-connected subgraph which embeds in $Q_{n-2}$ whenever $n \ge 4$.
\end{Theorem}
\begin{proof}
The graph is clearly a subgraph of $Q_{n-2}$ because all the words begin with $1$ and are of length $n-1$ so they lie in the subgraph of $Q_{n-1}$ isomorphic to $Q_{n-2}$. We now show that every such graph is $2$-connected. This is true for $n=4$. Let us assume this to be true for $n-1$, $n \ge 4$. Thus every word ending with $0$ obtained from $\Pi_{n-1} / \lbrace 0^{n-2}\rbrace$ will also be $2$-connected. Now consider any word ending in $1$, we already know it is connected to a word ending in $0$ and this is 2-connected to every other vertex with $0$. If the word ending in $1$ contains one or more $0$s then switch the leftmost $0$ in the string to a $1$. This will be an allowed word ending in $1$, if we switch the rightmost $1$ to a $0$ this will be a new vertex that is 2-connected to every other vertex ending in $0$, hence the chosen word is $2$-connected. If the word chosen consists only of $1$s, switch the second to last $1$ to a $0$. For any $n \ge 4$ this will be an allowed word ending in $1$, hence we are done.
\end{proof}

This result is the best possible because the vertex $1010^m$ for $m \ge 0$ is always at most $2$-connected in the partition graph.\par

\begin{Theorem}
The graph $\Pi_n / \lbrace 0^{n-1}\rbrace$ (i.e. with the $0$ vertex removed) is not hamiltonian for $n\ge 8$ but ${\Pi_n / \lbrace 0^{n-1}\rbrace}^2$ is hamiltonian for $n \ge 4$ hence the partitions words admit a $2$-Gray code.
\end{Theorem}
\begin{proof}  
The elements in each section of the bipartite division of the graph consist of partitions with an even or an odd number of parts. The difference between the two is given by the generating function, see OEIS A081362 $$\prod_{n=1}^{\infty}\frac{1}{1+q^n}.$$ The difference is $ > 1$ for $n > 7$. Since this is a bipartite graph, no hamiltonian cycle is possible. Now by Fleischner's Theorem ~\cite{Fleischner}, since ${\Pi_n / \lbrace 0^{n-1}\rbrace}$ is 2-connected for $n \ge 4$, the square graph ${\Pi_n / \lbrace 0^{n-1}\rbrace}^2$, formed by considering edges between all vertices of  ${\Pi_n / \lbrace 0^{n-1}\rbrace}$ separated by a Hamming distance of at most $2$, is hamiltonian.
\end{proof}

Herbert Wilf at the SIAM Conference on Discrete Math in June 1988 ~\cite{survey} had asked whether integer partitions of $n$ admit a Gray code. An allowed transitions would switch two parts $\lambda_i, \lambda_j$ of a partition of $n$ to $\lambda_i-1, \lambda_j+1$ as long as this switch also resulted in a valid partition of $n$. In addition, transitions that would increase or decrease a part of a partition by one and correspondingly remove or add a part of size one to this partition were also allowed. Savage ~\cite{savage} showed this was possible by working out an explicit algorithm. The 2-Gray code whose existence is proved in this paper involves transitions of only one kind, the collapse of a plus sign in a partition of $n$ or the insertion of a plus sign where permitted.\par

\section{Computational study of the edge distribution}\label{sec:compute}

Now consider the edges of the sequence of graphs. Since the graphs are subgraphs of the hypercube, the maximum possible number of edges for $\Pi_n / \lbrace 0^{n-1}\rbrace$ is bounded by $p(n)log(p(n))$ ~\cite{chung}. Two vertices are joined by an edge if they differ only in one place, i.e. if one of letters in a vertex at this slot is a $0$ and the other is $1$. Now consider each vertex consisting of at least one zero. Moving from left to right if we replace the first zero by a one (this is always possible) we get another vertex joined t0 this by an edge. Since the only vertex that does not have a zero is the all-ones vertex, we know that there are at least $p(n)-1$ edges. A subsequent zero in any vertex can be replaced by a one only if the numbers of ones before the zero preceding this zero is greater than the sum of the number of ones between the two zeros and the number of ones after the second zero. This can be reinterpreted in the conventional way of writing partitions, $\sum_{i=1}^{n}\lambda_i = n$, as the statement that there is one edge for every occurrence of $\lambda_i \geq \lambda_{i+1} + \lambda_{i+2}$ counted over all partitions on $n$. Determining the asymptotic value of this quantity remains an open question, and the graph theoretic result gives an upper bound of $\mathcal{O}(n^{1/2}p(n))$. \par

It is natural to try and understand the degree distribution of this network for large $n$.\par
\textbf{The bar graph below displays the degree sequence for $n=37$ i.e a graph consisting of $p(37) = 21637$ vertices (the all $0$s vertex is the only one with degree one).}

\resizebox{\textwidth}{!}{\includegraphics{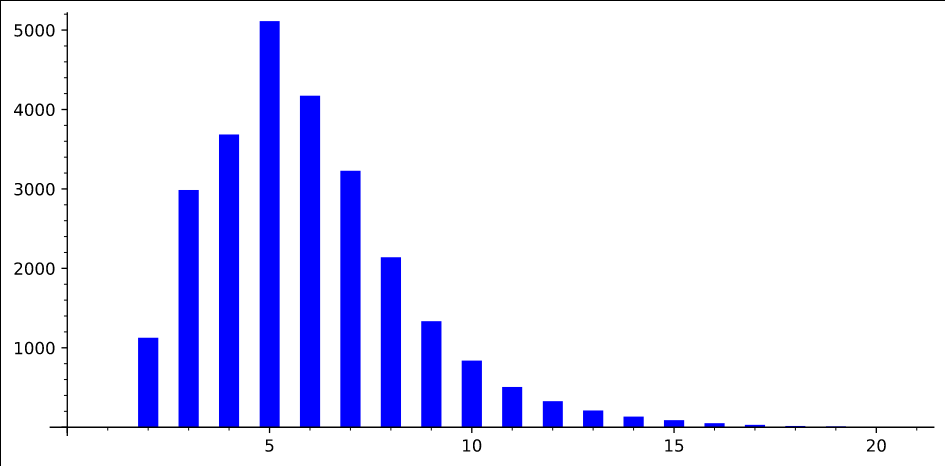}}

The evidence leads us to conjecture that the degree distribution of the partition graph tends to the lognormal distribution. This is supported by a number of different statistical tests on the data. The emergence of a longtailed distribution, often observed in many large networks, in this completely deterministic setting seems to be some interest. Given the nature of the growth of the graph sequence, it makes sense to study whether the degree distribution bears any relation to distribution of the number of zeros in the partition words of $n$. The evidence suggest this is indeed the case,  while the mean and variance differ, the distribution also seems to be a longtailed distribution best fitted by the lognormal distribution.\par

\textbf{The bar graph below illustrates the distribution of the number of partition words of $p(37)$ with a given number of zeros.} \par\noindent

\resizebox{\textwidth}{!}{\includegraphics{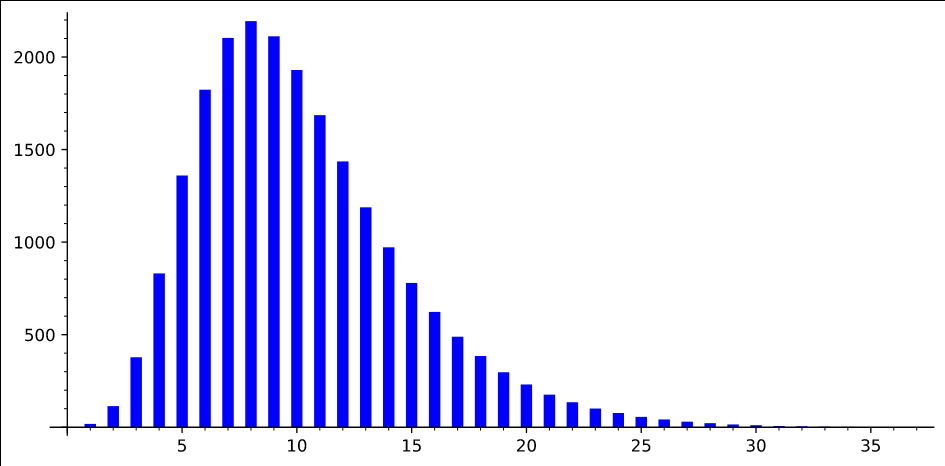}}

The possible distribution of the number of zeros of these partitions graphs leads to an interesting conjecture that goes back to the early history of the study of partitions initiated by Erdos ~\cite{erdos}, the distribution of the parts of a partitions. Considerable work has been done on this, with separate results for large parts of a partitions of the order of $n^{1/2}$ or small parts that tend to zero with $n$ ~\cite{fristedt} but there has been no description of a possible overall distribution. Since each zero corresponds to a plus sign in the canonical way of writing partitions, the number of zeros in a word is just one less than the number of parts in a graphs, hence both follow the same distribution. The evidence here suggests that the number of parts of a partition of $n$ tends to the lognormal distribution as $n$ becomes large.

\section{$d$-Dimensional Partitions as Words}

We first state the rules for plane partitions as they then extend easily to rules for all $d$-dimensional partitions. We then use the general case to show that these rules give a bijection between $d$-dimensional partitions and allowed words.\par  

\begin{Definition} \label{rule2}
The rules for representing plane partitions with ternary words on the alphabet $\lbrace2,1,0\rbrace$:\\
  a) Consecutive $0$s can only terminate a string.\\
  b) Words consisting of any two letters will follow the rules given for a binary alphabet with the understanding that the higher valued symbol will take the role of $1$ and the lower of $0$.\\
  c) Allowed words on the two symbols $\lbrace2,1\rbrace$ will be of the form\\
   $A=A_11A_21 \cdots A_i1 \cdots A_n1^m$, where the $A_i$ are a string of $2$s. Given two such words $B=B_11B_21 \cdots 1B_i1 \cdots 1B_n1^m$ and $C= C_11C_21 \cdots 1C_i1 \cdots 1C_r1^s$, we say $B0C$ is a word if $n \geq r $, $m + n \geq s + r $ and $B_i$ dominates $C_i$ for each $i \le r$. In such a case we say $B$ dominates $C$.\\
  d) Given any allowed word of the form $A_10A_20 \cdots A_i0 \cdots A_n$, where $A_i$ are allowed words on the alphabet $\lbrace2,1\rbrace$ , then $A_10A_20 \cdots A_i0 \cdots A_n0A_{n+1}$ is an allowed word if $A_n0A_{n+1}$ is an allowed word.\\
      \end{Definition}
We note that the definition implies that words such as $2210211$ are not allowed as the number of $1$s after the $0$ exceed the number of $1$s in the string preceding the $0$. From the rules suggest it is clear that a generic $2$-dimensional partition word will be of the form $$2^{a_{11}}12^{a_{12}}1 \cdots 12^{a_{1n}}1^{b_1}0   2^{a_{21}}12^{a_{22}}1 \cdots 12^{a_{2n}}1^{b_2}0 \cdots 02^{a_{j1}}12^{a_{j2}}1 \cdots 12^{a_{jn}}1^{b_j}0^m,$$ where $a_{ij} \ge 0$,  $a_{ij} \geq a_{(i+1)j}$, $a_{ij} \geq a_{i(j+1)}$ for all $i,j$ and $b_i \ge b_{i+1}$ for all $i$ .\par 

We recursively enumerate the plane partitions for $n \le 5$. Beginning with the empty word we obtain all the plane partitions of $n$ by appending $0$ to all the plane partitions of $n-1$, and then appending $1$ or $2$ as the rules listed above allow. Given the plane partition words for $n=3$ - $00, 10, 20, 11, 21, 22$, we append $0$ to each giving $000$, $100$, $200$, $110$, $210$, $220$, $101$ and then we append $1$, giving $101$, $111$, $211$, $221$ as words of the form $001$ or $201$ are not allowed and finally we append $2$ giving us $202$, $212$, $222$ as words of the form $002$, $102$ or $112$ are not allowed.\\

$$
\begin{array}{| c | l | c |}\hline
n & \text{Allowed words of length  } n-1 & p_2(n) \\ \hline\hline
1 & \epsilon & 1\\
2 &  0, 1, 2 & 3 \\
3 & 00, 10, 20, 11, 21, 22 & 6 \\
4 &  \begin{minipage}[t]{0.75\textwidth}
$000$, $100$, $200$, $110$, $210$, $220$, $101$, $111$, $211$, $221$, $202$, $212$, $222$ 
\end{minipage} & 13\\
5 &  \begin{minipage}[t]{0.75\textwidth}
$0000$, $1000$, $2000$, $1100$, $2100$, $2200$, $1010$, $1110$, $2110$, $2210$, $2020$, $2120$, $2220$, $1101$, $2101$, $1111$, $2111$, $2211$, $2121$, $2221$, $2102$, $2202$, $2212$, $2222$
\end{minipage} & 24\\
\hline
\end{array}
 $$ \\

We now write out the rules for $d$-dimensional partition words assuming we have satisfactorily defined $d-1$-dimensional partition words.\\
   \begin{Definition}\label{rule3}
The rules for representing $d$-dimensional partitions, for $d \ge 3$, with words on the alphabet $\lbrace d, d-1, \cdots , 2, 1, 0\rbrace$:\\
  a) Consecutive $0$s can only terminate a string.\\
  b) Words consisting of any $d$ letters  will follow the rules for $(d-1)$-dimensional partition words with the understanding that the symbols will take the role of\\
   $\lbrace d-1, \cdots , 2, 1, 0 \rbrace$ in $(d-1)$-dimensional partition words in decreasing order of their numerical value.\\
  c) Allowed words on the $d$ symbols $\lbrace d, d-1, \cdots , 2, 1\rbrace$ will be of the form\\
   $A=A_11A_21 \cdots A_i1 \cdots A_n1^m$, where the $A_i$ are a string on the $d-2$ symbols $\lbrace d, d-1, \cdots , 2\rbrace$. Given two such words $B=B_11B_21 \cdots 1B_i1 \cdots 1B_n1^m$ and $C= C_11C_21 \cdots 1C_i1 \cdots 1C_r1^s$, we say $B0C$ is a word if $n \geq r $, $m + n \geq s + r $ and $B_i$ dominates $C_i$ for each $i \le r$. In such a case we say $B$ dominates $C$.\\
  d) Given any allowed word of the form $A_10A_20 \cdots A_i0 \cdots A_n$, where $A_i$ are allowed words on the alphabet $\lbrace d, d-1, \cdots ,2,1 \rbrace$ , $A_10A_20 \cdots A_i0 \cdots A_n0A_{n+1}$ is an allowed word if $A_n0A_{n+1}$ is an allowed word.\\
    \end{Definition}\par

\begin{Theorem}\label{all}  The rules listed above ~\eqref{rule3} give a bijection between $d$-dimensional partitions of $n$ and allowed words of length $n-1$ on an alphabet of $d+1$ letters.
\end{Theorem}
\begin{proof}
We have shown the result if true for $n=1$ and $n=2$. We assume it is true for $d=m-1$, and we show it is true for $d=m$. An $m$-dimensional partition of $n$ is given by $n_{i_1, \cdots ,i_m}$ with $\sum_{i_1 \leq r_1, \cdots , i_m \leq r_m} n_{i_1, \cdots, i_m} = n$ such that $n_{i_1, \cdots ,i_m} \geq n_{j_1, \cdots, j_m}$ for $i_1 \leq j_1, \cdots , i_m \leq j_m$.\par 
By assumption we have the fact that for each $i_m$ we can write the partition $\sum_{i_1, \cdots ,i_{m-1}} n_{i_1, \cdots, i_m} = n_{i_m}$ uniquely as a partition word $A_{i_m}$ on the $m$ letter alphabet $\lbrace m, m-1, \cdots ,1\rbrace$. We assume without loss of generality that for some $1 \leq t \leq r_m$, all the $n_{k_m}$ for $k_m \geq t$ are $1$. Hence the $m$-dimensional partition word $$A_10A_20A_3 \cdots A_{t-1}0^{r_m-t}$$ corresponds uniquely to the given partition.
\end{proof}

We illustrate the ease of enumeration of higher-dimensional partitions by working out the $3$-dimensional partitions up to $n=6$. This simple calculation carried out by hand is of some interest as MacMahon's conjectured a generating function for such partitions that was first shown to fail for $n=6$ in \cite{atkin}. The book by Andrews \cite{andrews} remains the standard reference for this calculation, and in fact the study of solid partitions.\par

Beginning with the empty word we obtain all the solid partitions of $n$ by appending $0$ to all the solid partitions of $n-1$, and then appending $1$, $2$ or $3$ as the rules listed above allow. We note that the rules do not allow words of the form $3102$ or $3320312$.\par

$$
\begin{array}{| c | l | c |}\hline
n & \text{Allowed words of length  } n-1 & p_3(n) \\ \hline\hline
1 & \epsilon & 1\\ 
2 &  0, 1, 2, 3 & 4 \\ 
3 & 00, 10, 20, 30, 11, 21, 31, 22, 32, 33  & 10 \\
4 &  \begin{minipage}[t]{0.75\textwidth}
$000$, $100$, $200$, $300$, $110$, $210$, $310$, $220$, $320$, $330$, $101$, $111$, $211$, $311$, $221$, $321$, $331$, $202$, $212$, $222$, $322$, $332$, $303$, $313$, $323$, $333$ 
\end{minipage} & 26\\
5 &  \begin{minipage}[t]{0.75\textwidth}
$0000$, $1000$, $2000$, $3000$, $1100$, $2100$, $3100$, $2200$, $3200$, $3300$, $1010$, $1110$, $2110$, $3110$, $2210$, $3210$, $3310$, $2020$, $2120$, $2220$, $3220$, $3320$, $3030$, $3130$, $3230$, $3330$, $1101$, $2101$, $3101$, $1111$, $2111$, $3111$, $2211$, $3211$, $3311$, $2121$, $2221$, $3221$, $3321$, $3131$, $3231$, $3331$, $2102$, $2202$, $3202$, $2212$, $3212$, $2222$, $3222$, $3322$, $3232$, $3332$, $3103$, $3203$, $3303$, $3213$, $3313$, $3323$, $3333$ 
\end{minipage} & 59\\
6 &  \begin{minipage}[t]{0.75\textwidth}
$00000$, $10000$, $20000$, $30000$, $11000$, $21000$, $31000$, $22000$, $32000$, $33000$, $10100$, $11100$, $21100$, $31100$, $22100$, $32100$, $33100$, $20200$, $21200$, $22200$, $32200$, $33200$, $30300$, $31300$, $32300$, $33300$, $11010$, $21010$, $31010$, $11110$, $21110$, $31110$, $22110$, $32110$, $33110$, $21210$, $22210$, $32210$, $33210$, $31310$, $32310$, $33310$, $21020$, $22020$, $32020$, $22120$, $32120$, $22220$, $32220$, $33220$, $32320$, $33320$, $31030$, $32030$, $33030$, $32130$, $33130$, $33230$, $33330$, $10101$, $11101$, $21101$, $31101$, $22101$, $32101$, $33101$, $21201$, $31301$, $11011$, $11111$, $21111$, $31111$, $22111$, $32111$, $33111$, $21211$, $22211$, $32211$, $33211$, $31311$, $32311$, $33311$, $21021$, $22121$, $32121$, $22221$, $32221$, $33221$, $32321$, $33321$, $31031$, $32131$, $33131$, $33231$, $33331$, $21102$, $22102$, $32102$, $20202$, $21202$, $22202$, $32202$, $33202$, $32302$, $21212$, $22212$, $32212$, $33212$, $32312$, $22022$, $22122$, $22222$, $32222$, $33222$, $32322$, $33322$, $32032$, $32132$, $33232$, $33332$, $31103$, $32103$, $33130$, $32203$, $33203$, $30303$, $31303$, $32303$, $33303$, $32213$, $33213$, $31313$, $32313$, $33313$, $32323$, $33323$, $33033$, $33133$, $33233$, $33333$
\end{minipage} & 140\\
\hline
\end{array}
 $$

\section{d-Partition Subgraphs of Hamming Graphs}\label{hamming}

\begin{Theorem}
$\Pi(d,n)$ is a connected (d+1)-partite subgraph of the hamming graph $\textit{H}{(n-1,d+1)}$. 
\end{Theorem}
\begin{proof}
The d-partition graph for $n=2$, consisting of words of length $1$, is the complete graph $\textit{K}_{d+1}$. By the recursive definition for generating these d-partitions, it is clear that $\Pi(d,i)$ can be seen as  a subgraph of $\Pi(d,i+1)$ by appending $0$ to each of its vertices, so $\textit{K}_{d+1}$ is a subgraph of $\Pi(d,i)$ for all $i$. Hence, it requires at least $d+1$ colours so that no adjacent vertices are assigned the same colour. Moreover, since each $\Pi(d,i)$ is a subgraph of  $\textit{H}{(i-1,d+1)}$ by definition, it inherits the $d+1$ coloouring of the hamming graph, hence it is $d+1$-partitite.\par

Now, If we assume the result is true for $n-1$, i.e. words of length $n-2$, then the partition words of length $n-1$ are of two types, those obtained from partition words of length $n-2$ by appending $0$ to each word, and those obtained by appending $1$ where allowed. The words ending in $0$ clearly form a connected subgraph of $Q_{n-1}$. Moreover, any word ending in $1$ will also have an associated word ending in $0$ which only differs at the last place, and hence it will be connected to the subgraph by a vertex of the $Q_{n-1}$ hypercube. It is also clear that the graph is bipartite because the words with an odd(even) number of $1$s only connect the words with an even(odd) number of $1$s.
\end{proof}

\begin{Theorem}
The diameter of $\Pi(d,n)$ is at most $2n-3$.
\end{Theorem} 
\begin {proof}
Let $A$ and $B$ be any two vertices of $\Pi(d,n)$. Moving right to left along the word $A$, switch each letter to a $0$, it is easy to see from the rules that the new word so formed at each step is a valid d-partition word. Continue this process till only the starting letter is non-zero. Flip the starting letter to the starting letter of $B$, and moving left to right along this word switch the zero to the letter in the corresponding position in $B$. Again, the resulting word is a d-partition word at each step. The total number of flips requires is one less than the the total length of both the word, $(2(n-1) + 1$.\par

The distance between the vertices $A=0^{n-1}$ and $B=1^{n-1}$ is $n-1$. We note that any vertex has a connected path to $A$ through allowed partition words by moving right to left in the string and changing each $1$ to a $0$ or to $B$ by moving left to right in the string and changing each $0$ to a $1$. Let $x_1$ and $x_2$ be any two vertices of the graph. If the distance of $x_1$ from $A$ is $a_1$, then the distance of $x_1$ from $B$ is $n-1-a_1$, and one of the two terms is $\le$ $\lfloor{(n-1)/2}\rfloor$, similarly the distance of $x_2$ from $A$ is $b_1$, then the distance of $x_2$ from $B$ is $n-1-a_2$, and one of them is $\le$  $\lfloor{(n-1)/2}\rfloor$. By choosing the appropriate segment joining $x_1$ to $A$ or $B$, as the case may be, it follows that the shortest path between $x_1$ and $x_2$ via $A$ or $B$ is $\leq$ ${n-1}$.
\end{proof}

\begin{Theorem}
The graph $\Pi(d,n)$ is connected for all $d$ and $n$.
\end{Theorem}
\begin{proof}
We have shown the graph $\Pi(1,n)$ is connected for all $n$ (the $0^{n-1}$ partition is connected to $10^{n-2}$). Now consider the case of fixed $d$ as $n$ varies. For $n=2$ the graph is connected. We proceed by induction and assume the result is true for fixed $d$ and all $n \le M$, now consider the graph formed by all the partitions of $\Pi(d,M+1)$. All the partition words ending with $0$ form the same graph as $\Pi(d,M)$ hence it is connected. Consider any partition ending with any of $d$ letters $\ne 0$, then this partition word is connected to the word which only differs in the rightmost place where it has a $0$, hence it is connected to a connected graph so the new graph is also connected.
\end{proof}

\begin{Theorem}
The graph $\Pi(d,n)$ is not 2-connected for all $d\ge2$ and $n\ge6$.
\end{Theorem}
\begin{proof}
For $d=2$ and $n=6$, the partition word $21021$ is only connected to $21020$. We can append a $0$ to it and show this is true for $n=7$. For $n=8$ , we consider $2210221$ and $2110211$. The process can be thus extended indefinitely. Moreover, these partition words are valid partition words for all $d\ge2$ and $n\ge6$, hence the result follows. Unlike the case for ordinary partitions, these vertices connected by only one edge to graphs increase in number as $n$ increases and we cannot think of the graph as $2$-connected just by isolating a handful of vertices.
\end{proof}

\textbf{The graph for $\Pi(2,5)$, consisting of 24 vertices, is shown below, and is the largest 2-connected graph in the sequence $\Pi(2,n)$.}\par

\resizebox{\textwidth}{!}{\includegraphics{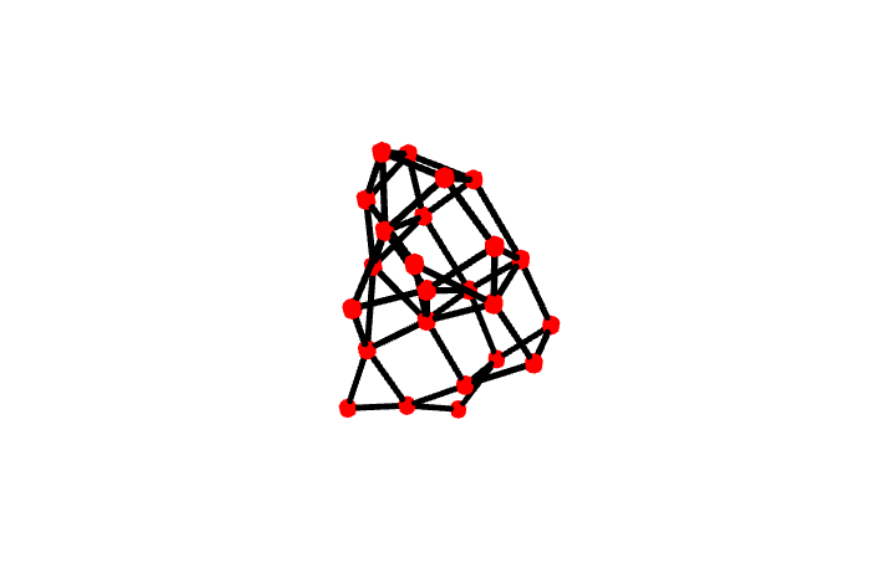}}

By an old result of Sekanina ~\cite{sekina} it follows that the cube of a connected graph admits a Hamiltonian cycle hence this implies:

\begin{Theorem}

d-dimensional Partitions admit a 3-Gray code.

\end{Theorem}

These graphs also tend to show the same lognormal distribution for large $n$. Unlike the binary case, which yields bipartite graphs with low clustering coefficients, the clustering coefficient seems to increase with $d$, suggesting the possibility of modelling large networks using variations of such a deterministic scheme.\par

\section{Acknowledgements}

I would like to thank Gaurav Bhatnagar for discussions and insights into much of the work in this paper, as well as help with programming. All numerical computations and graphics have been generated using Sage ~\cite{sage}.

\end{document}